\documentclass[12pt]{article}

\usepackage{amsmath}\usepackage{amssymb}

\newtheorem{theorem}{Theorem}[section]
\newtheorem{lemma}{Lemma}[section]

\newtheorem{remark}{Remark}[section]
\newtheorem{example}{Example}[section]
\newtheorem{definition}{Definition}[section]

\newcommand{\Pol}{\mathop{\mathrm{Pol}}}
\newcommand{\Int}{\mathop{\mathrm{Int}}}
\newcommand{\B}{\mathrm{B}}
\newcommand{\C}{\mathrm{C}}

\renewcommand{\o}{\mathrm{o}}

\newcommand{\F}{\mathrm{F}}
\newcommand{\SQ}{\mathrm{SQ}}

\newcommand{\N}{\mathbb{N}}
\newcommand{\R}{\mathbb{R}}
\newcommand{\Z}{\mathbb{Z}}
\newcommand{\D}{\Delta}

\title{Regional topology and approximative solutions 
of difference and differential equations} 
\author{Janusz Migda}
 
\begin{document}\maketitle

\begin{abstract}
We introduce a topology, which we call the regional topology, 
on the space of all real functions on a given locally compact metric space. Next we 
obtain a new versions of Schauder's fixed point  theorem and Ascoli's theorem. We use 
these theorems and the properties of the iterated remainder operator to establish  conditions under which there exist solutions, with prescribed asymptotic behavior, of 
some difference and differential equations. 
\medskip\\
{\bf Key words:} regional topology, difference equation, differential equation, 
remainder operator, approximative solution, prescribed asymptotic behavior, 
Schauder's theorem, Ascoli's theorem.\\
{\bf AMS Subject Classification:} $39A10$
\end{abstract}

\section{Introduction}

Let $\N$, $\Z$, $\R$ denote the set of positive integers, the set of all integers 
and the set of real numbers, respectively. Assume 
\[
m\in\N, \quad t_0\in[0,\infty), \quad I=[t_0,\infty).
\]
In this paper we consider difference equations of the form 
\begin{equation}\label{E1}\tag{E1}
\D^mx_n=a_nf(n,x_n)+b_n 
\end{equation}
\[
n\in\N,\quad a_n,\ b_n \in\R,\quad f:\N\times\R\to\R,
\]
and differential equations of the form 
\begin{equation}\label{E2}\tag{E2}
x^{(m)}(t)=a(t)f(t,x(t))+b(t) 
\end{equation}
\[
t\in I,\quad a(t),\ b(t)\in\R,\quad f:I\times\R\to\R.
\]
We are interesting in the existence of solutions with prescribed asymptotic behavior. 
More precisely, we establish conditions under which, for a given $\alpha\in(-\infty,0]$ 
and solution $y$ of the equation $\D^my=b$, there exists a solution $x$ of \eqref{E1} 
such that $x=y+\o(n^\alpha)$. Analogously, we establish conditions under which, for 
a given $\alpha\in(-\infty,0]$ and solution $y$ of the equation $y^{(m)}=b$ there 
exists a solution $x$ of \eqref{E2} such that $x=y+\o(t^\alpha)$. 

We obtain these results using a new version of the Schauder fixed point theorem and 
a new version of the Ascoli theorem. First, we introduce the notion of regional norm 
and regional topology. In the case of functional space, the regional topology is the 
topology of uniform convergence. This topology is nonlinear, but it is almost 
linear, see Remark \ref{top 03}. Using this topology we can state our version of the  Schauder theorem. 

Next, in functional spaces, we introduce the notion of homogeneous at infinity and 
stable at infinity family of functions. The first notion generalizes the notion of 
equiconvergence at infinity, see Remark \ref{homogeneous} and Remark 
\ref{equiconvergence}. We use the notion of homogeneity at infinity to state our 
version of the Ascoli theorem. 

Our approach to the study of existence of solutions with prescribed asymptotic 
behavior is based on using the iterated remainder operator. In the discrete case 
properties of this operator were presented in \cite{J Migda 2010 a}, 
\cite{J Migda 2014 a} and \cite{J Migda 2014 c}. The basic properties of this 
operator in the continuous case we establish in the first part of Section 5. 

If $y$ is a sequence with known properties and $x$ is a solution of \eqref{E1} such 
that $x=y+\o(1)$, then we say that $x$ is a solution with prescribed asymptotic 
behavior. Writing the equality $x=y+\o(1)$ in the form $y=x+\o(1)$ we may say 
that $y$ is an approximative solution of \eqref{E1}. 

In Section 4, using the iterated remainder operator we establish conditions under 
which a given solution $y$ of the equation $\D^my=b$ is an approximative solution of  \eqref{E1}. Moreover, the technique of iterated remainder operator allow us to change  $\o(1)$ by $\o(n^\alpha)$ for a given nonpositive real $\alpha$. Hence we obtain an 
approximative solution which better approximates the solution $x$. 

In Section 5, we obtain an analogous result in the continuous case. In Section 6, 
we give some additional remarks on the regional topology. 

This paper is a continuation of a cycle of papers 
\cite{MJ Migda 2001}-\cite{J Migda 2014 b}. Our studies, in continuous case, were 
inspired by papers \cite{Ehrnstrom, Lipovan, Mingarelli, Mustafa 2006} 
(devoted to asymptotically linear solutions) and by 
papers \cite{Hallam, Kong, Philos 2004, Trench 1984} 
(devoted to asymptotically polynomial solutions).

\section{Notation and terminology}

If $p,k\in\Z$, $p\leq k$, then $\N_p$, $\N(p,k)$ denote the sets defined by 
\[
\N_p=\{p,p+1,\dots\}, \qquad  \N(p,k)=\{p,p+1,\dots,k\}.
\]
Let $X$ be a set. We denote by 
\[
\F(X), \quad \F_b(X)
\]
the space of all functions $f:X\to\R$ and the space of all bounded functions $f\in\F(X)$ 
respectively. Let $f,g\in\F(X)$, $F\subset\F(X)$. Then 
\[
|f|, \quad f+g, \quad f-g, \quad fg
\]
denotes the functions defined by a standard way. If $Z\subset X$, then $f|Z$ denotes 
the restriction $f|Z:Z\to\R$. Moreover, 
\[
F-F=\{f-g:\ f,g\in F\}, \qquad F+g=\{f+g:\ f\in F\}.
\]
We say that $F$ is pointwise bounded if for any $t\in X$ the set 
$F(t)=\{f(t):\ f\in F\}$ is bounded. If $X$ is a topological space, then $\C(X)$ 
denotes the space of all continuous functions $f\in\F(X)$.
\smallskip\\
Assume $X$ is a locally compact metric space. A function $f\in\F(X)$ is called 
vanishing at infinity if for any $\varepsilon>0$ there exists a compact subset $Z$ of 
$X$ such that $|f(t)|<\varepsilon$ for any $t\notin Z$. The space of all vanishing at  infinity functions $f\in\F(X)$ we denote by $\F_0(X)$. Moreover 
\[
\C_0(X)=\C(X)\cap\F_0(X). 
\]
Note that any continuous and vanishing at infinity function is bounded. 
\smallskip\\
Assume $p\in\N$, $t_0\in[0,\infty)$, $I=[t_0,\infty)$. Then 
\[
\F_0(\N_p)=\{x\in\F(\N_p):\ \lim_{n\to\infty}x_n=0\}, \quad 
\F_0(I)=\{f\in\F(I):\ \lim_{t\to\infty}f(t)=0\}.
\]
Moreover, for $\alpha\in\R$ we use the following notation 
\[
\o(n^\alpha)=\{x\in\F(\N_p):\ \lim_{n\to\infty}x_nn^{-\alpha}=0\}, \qquad 
\o(t^\alpha)=\{f\in\C(I):\ \lim_{t\to\infty}f(t)t^{-\alpha}=0\}
\]
Let $X$ be a metric space. For a subset $A$ of $X$ and $\varepsilon>0$, we define 
an $\varepsilon$-framed interior of $A$ by 
\[
\Int(A,\varepsilon)=\{x\in X:\ \mathrm{B}(x,\varepsilon)\subset A\},
\]
where $\B(x,\varepsilon)$ denotes an open ball of radius $\varepsilon$ 
about $x$. Moreover, we define an $\varepsilon$-ball about $A$ by 
\[
\B(A,\varepsilon)=\bigcup_{x\in A}\B(x,\varepsilon). 
\]
A subset $A$ of $X$ is called an $\varepsilon$-net for a subset $Z$ of $X$ if 
$Z\subset\B(A,\varepsilon)$. A subset $Z$ of $X$ is said to be totally bounded 
if for any $\varepsilon>0$ there exist a finite $\varepsilon$-net for $Z$. 
\smallskip\\
Let $d$ denotes a metric in $X$. A family $F\subset\F(X)$ is called equicontinuous 
if for any $\varepsilon>0$ there exists a $\delta>0$ such that the condition $s,t\in X$,  $d(s,t)<\delta$ implies $|f(s)-f(t)|<\varepsilon$ for any $f\in F$. 
For $m\in\N$ we use the rising factorial notation 
\[
n^{\overline{m}}=n(n+1)\dots(n+m-1)\quad \text{with} \quad n^{\overline{0}}=1.
\]

\section{Regional topology}

Let $X$ be a set. For a function $f\in\F(X)$ we define a generalized norm  $\|f\|\in[0,\infty]$ by 
\[
\|f\|=\sup\{|f(t)|:\ t\in X\}.
\]
We say that a subset $F$ of $\F(X)$ is ordinary if $\|f-g\|<\infty$ 
for any $f,g\in F$. We regard every ordinary subset $F$ of $\F(X)$ 
as a metric space with metric defined by 
\begin{equation}\label{dfg}
d(f,g)=\|f-g\|.
\end{equation}
Let $U\subset\F(X)$. We say that $U$ is regionally open if $U\cap F$ 
is open in $F$ for any ordinary subset $F$ of $\F(X)$. The family of 
all regionally open subsets is a topology on $\F(X)$ which we call the regional 
topology. We regard any subset of $F(X)$ as a topological space with topology induced 
by the regional topology. 
\smallskip\\
Note that a subset $U$ of $\F(X)$ is a neighborhood of $f\in U$ if and only if there 
exists an $\varepsilon>0$ such that $\B(f,\varepsilon)\subset U$. Hence, a sequence 
$(f_n)$ in $\F(X)$ is convergent to $f$ if and only if for any the sequence 
$\|f_n-f\|$ is convergent to zero. Therefore we can say that the regional 
topology in $\F(X)$ is the topology of uniform convergence. 
\smallskip\\
More generally, let $X$ be a real vector space. We say that a function 
$\|\cdot\|:X\to[0,\infty]$ is regional norm if the condition $\|x\|=0$ is 
equivalent to $x=0$ and for any $x,y\in X$ and $\alpha\in\R$ we have 
\[
\|\alpha x\|=|\alpha|\|x\|, \quad \|x+y\|\leq\|x\|+\|y\|.
\]
Hence, the notion of regional norm generalizes the notion of usual norm. 
If a regional norm on $X$ is given, then we say that $X$ is a regional 
normed space. If there exists a vector $x\in X$ such that $\|x\|=\infty$, then we 
say that $X$ is extraordinary. 
\smallskip\\
Assume $X$ is a regional normed space. 
We say that a subset $Z$ of $X$ is ordinary if $\|x-y\|<\infty$ 
for any $x,y\in Z$. We regard every ordinary subset $Z$ of $X$ 
as a metric space with metric defined by \eqref{dfg}. Analogously as above we
define a regional topology on $X$. Let 
\[
X_0=\{x\in X:\ \|x\|<\infty\}.
\]
Obviously $X_0$ is a linear subspace of $X$. Moreover, the regional norm 
induces a usual norm on $X_0$. We say that $X$ is a Banach regional 
space if $X_0$ is complete.
\smallskip\\
For $p\in X$ let 
\[
X_p=p+X_0.
\]
Then $X_p$ is a maximal ordinary subset of $X$ which contain $p$. It is easy to see 
that $X_p$ is the equivalence class of $p$ under the relation defined by 
\[
x\equiv y \Leftrightarrow \|x-y\|<\infty.
\]
We say that $X_p$ is an ordinary component or a region of the space $X$. It is easy 
to see that $X_p$ is also a connected component of $p$ in $X$. From topological point 
of view, the space $X$ is a disjoint union of all its regions. In particular, every 
region $X_p$ is open and closed subset of $X$. Moreover, for any $p\in X$ the 
translation 
\[
T_p:X_0\to X_p, \qquad T_p(x)=x+p
\]
is an isometry of $X_0$ onto $X_p$. Hence a metric space $X_p$ is metrically equivalent 
to the normed space $X_0$. Note also, that the translation $T_p$ preserves convexity 
of subsets.

Now, we are ready to state and prove a generalization of the following theorem.
\medskip\\
\textbf{Theorem}\ \textbf{(Schauder fixed point theorem)} 
\textit{Assume $Q$ is a closed and convex subset of a Banach space $X$, 
a map $A: Q\to Q$ is continuous and the set $AQ$ is totally 
bounded. Then there exists a point $x\in Q$ such that $Ax=x$.} 

\begin{theorem}[Generalized Schauder theorem]\label{GST}
Assume $Q$ is a closed and convex subset of a regional Banach space $X$, 
a map $A: Q\to Q$ is continuous and the set $AQ$ is ordinary and totally 
bounded. Then there exists a point $x\in Q$ such that $Ax=x$. 
\end{theorem}
\textbf{Proof.} 
Choose $p\in AQ$. Since $AQ$ is ordinary, we have 
\begin{equation}\label{AQ}
AQ\subset X_p. 
\end{equation}
Let 
\[
Q_p=Q\cap X_p, \qquad Q_0=Q_p-p.
\]
Then $Q_p$ is closed and convex. Moreover, by \eqref{AQ}, $AQ_p\subset Q_p$. Let 
\[
A':Q_p\to Q_p, \qquad A'x=Ax, \qquad U:Q_p\to Q_0, \qquad Ux=x-p.
\]
Then $Q_0=T_p^{-1}Q_p$, and 
\[
T_p^{-1}: X_p\to X_0 
\]
is an isometry and preserves convexity. Hence $Q_0$ is a closed and convex subset 
of a Banach space $X_0$. Let 
\[
B:Q_0\to Q_0, \qquad B=U\circ A'\circ U^{-1}.
\]
Then $B$ is continuous and 
\begin{equation}\label{BQ}
BQ_0=(U\circ A'\circ U^{-1})Q_0=(U\circ A')(U^{-1}Q_0)=(U\circ A')Q_p=U(A'Q_p).
\end{equation}
Moreover, 
\[
A'Q_p=AQ_p\subset AQ.
\]
Hence $A'Q_p$ is totally bounded and $U$ is an isometry. Therefore, by \eqref{BQ}, 
$BQ_0$ is totally bounded. Thus, by the Schauder fixed point theorem, 
there exists a point $y\in Q_0$ such that $By=y$. Let $x=U^{-1}y$. Then 
\[
x=U^{-1}y=U^{-1}By=U^{-1}UA'U^{-1}y=A'U^{-1}y=A'x=Ax.
\]
The proof is complete. \quad $\Box$.

\begin{remark} 
In Example \ref{ord ess} we show, that the assumption of ordinarity of the set $AQ$, 
in the above theorem, is essential. 
\end{remark}
Note that, as in the functional case, a subset $U$ of $X$ is a neighborhood of a point 
$x$ if and only if there exists an $\varepsilon>0$ such that $\B(x,\varepsilon)\subset U$.  Hence, we have the following remark, which will be used in the proof of Theorem \ref{rho}.

\begin{remark}\label{continuity} 
If $T$ is a topological space, then a map $\varphi:X\to T$ is continuous at a 
point $x$ if and only if for any neighborhood $V$ of $\varphi(x)$ there 
exists an $\varepsilon>0$ such that $\varphi(\B(x,\varepsilon))\subset V$.
\end{remark}
Now, assume that $X$ is a locally compact metric space and $F\subset\F(X)$. 
We say that $F$ is: 
\begin{itemize}
\item[] vanishing at infinity if for any $\varepsilon>0$ there exists a compact                 subset $Z$ of $X$ such that $|f(s)|<\varepsilon$ for every                    $s\notin Z$ and every $f\in F$. 
\item[] equicontinuous at infinity if for any $\varepsilon>0$ there exists a compact            subset $Z$ of $X$ such that $|f(s)-f(t)|<\varepsilon$ for every               $s,t\notin Z$ and every $f\in F$. 
\item[] stable at infinity if for any $\varepsilon>0$ there exists a compact 
        subset $Z$ of $X$ such that $|f(s)-g(s)|<\varepsilon$  for any                $s\notin Z$ and $f,g\in F$
\item[] homogeneous at infinity if for any $\varepsilon>0$ there exists a compact               subset $Z$ of $X$ such that $|(f-g)(s)-(f-g)(t)|<\varepsilon$   
        for every $s,t\notin Z$ and every $f,g\in F$. 
\end{itemize}

\begin{remark}\label{homogeneous} 
Obviously, every family which is vanishing at infinity is also stable at infinity 
and equicontinuous at infinity. Moreover, every stable at infinity and every 
equicontinuous at infinity family is homogeneous at infinity. It is easy to see that 
a family $F\subset\F(X)$ is stable at infinity if and only if $F-F$ is vanishing at 
infinity. Similarly $F$ is homogeneous at infinity if and only if $F-F$ is 
equicontinuous at infinity. 
\end{remark}

\begin{example} 
Let 
\[
f,h:\R\to\R, \qquad f(t)=\exp(-t^2), \qquad h(t)=t^2,
\]
\[
F=\{f\}, \qquad G=\{f, f+1\}, \qquad H=\{h\}, \qquad K=\{h, h+1\}.
\]
Then $F$ is vanishing at infinity, $G$ is equicontinuous at infinity but not vanishing 
at infinity. $H$ is stable at infinity but not equicontinuous at infinity, $K$ is  homogeneous at infinity but not stable at infinity and not equicontinuous at infinity.
\end{example} 

In \cite{Avramescu 2002}, Avramescu use the term evanescent solution for solution 
vanishing at infinity. 
If $X=\N$, then the notion of equicontinuity at infinity is equivalent to the notion 
of uniformly Cauchy family of sequences, see \cite{Cheng 1993}.

Now, we are ready to generalize the following theorem.
\smallskip\\
\textbf{Theorem (Ascoli theorem)} 
\textit{If $X$ is a compact metric space, then every pointwise bounded and  
equicontinuous subset $F$ of $\C(X)$ is totally bounded. }

\begin{theorem}[Generalized Ascoli theorem]\label{GAT}
If $X$ is a locally compact metric space, then every equicontinuous, pointwise bounded  
and homogeneous at infinity subset $F$ of $\C(X)$ is totally bounded. 
\end{theorem}
\textbf{Proof.} 
Let $h\in F$ and 
\[
G=F-h=\{f-h: f\in F\}.
\]
Then $G$ is pointwise bounded, equicontinuous and equicontinuous at infinity.
Let $\varepsilon>0$. Choose a compact $Z\subset X$ such that 
\[
|g(s)-g(t)|<\varepsilon
\] 
for $g\in G$ and $s,t\notin Z$. 
Choose $s\in X\setminus Z$ and let 
\[
Y=Z\cup\{s\}. 
\]
Then $Y$ is a compact subset od the space $X$. By the Ascoli theorem, there exist 
\[
g_1,\dots,g_n\in G
\]
such that the family $g_1|Y,\dots,g_n|Y$ is an $\varepsilon$-net for the set 
$\{g|Y: g\in G\}$. Let $g\in G$. Then there exists an index $i\in\N(1,n)$ such that  $|g(y)-g_i(y)|\leq\varepsilon$ for any $y\in Y$. Let $t\in X\setminus Y$. Then 
\[
|g(t)-g_i(t)|\leq|g(t)-g(s)|+|g(s)-g_i(s)|+|g_i(s)-g_i(t)|\leq 3\varepsilon.
\]
Hence $\{g_1,\dots,g_n\}$ is a $3\varepsilon$-net for $G$. Therefore 
the set  $G$ is totally bounded. Thus the set $F=G+h$ is also 
totally bounded. \quad $\Box$. 

Theorem \ref{GAT} generalizes the Compactness criterion of C. Avramescu (see 
\cite[page 1164]{Philos 2004} and Remark \ref{equiconvergence}). In this paper, we 
do not need the Arzela theorem which states that if $X$ is a compact metric space, 
then every totally bounded subset of $\C(X)$ is equicontinuous. 

Now, we present the last `topological' theorem which will be used in the proofs of 
our main theorems.

\begin{theorem}\label{rho}
Let $X$ be a locally compact metric space, $h\in\F(X)$, $\rho\in\F_0(X)$ and 
\[
F=\{f\in\F(X):\ |f-h|\leq|\rho|\}.
\]
Then $F$ is closed, convex, pointwise bounded and stable at infinity. Moreover 
\begin{enumerate}
\item[$(a)$] if $\rho$ is continuous, then $F$ is ordinary, 
\item[$(b)$] if $X$ is discrete, then $F$ is compact and ordinary. 
\end{enumerate} 
\end{theorem}
\textbf{Proof.} 
Obviously, $F$ is convex, pointwise bounded and stable at infinity. 
Using Remark \ref{continuity}, it is easy to see that the map 
\[
\F(X)\to\F(X), \qquad f\mapsto f-h
\]
is continuous. Similarly, using Remark \ref{continuity}, we can see that 
for any $x\in X$, the evaluation 
\[
e_x:\F(X)\to\R, \qquad e_x(f)=f(x)
\]
is continuous. Hence, for any nonnegative $\alpha$, the set 
\[
\{f\in\F(X):\ |f(x)-h(x)|\leq\alpha\}
\]
is closed in $\F(X)$. Therefore the intersection 
\[
F=\bigcap\limits_{x\in X}\{f\in\F(X):\ |f(x)-h(x)|\leq|\rho(x)|\}
\] 
is closed. Assume $\rho$ is continuous, then $\|\rho\|<\infty$ and 
\[
\|f-g\|\leq\|f-h\|+\|g-h\|\leq 2\|\rho\|
\]
for any $f,g\in F$. Hence $F$ is ordinary.  
\smallskip\\ 
Now, assume that $X$ is discrete. 
Then $\rho$ is continuous and, by (a), $F$ is ordinary. Let $G$ be defined by 
\[
G=\{f\in\F(X):\ |f|\leq\rho\}.
\]
Choose an $\varepsilon>0$. Every compact subset of $X$ is finite. Since $\rho$ is 
vanishing at infinity, there exists a finite subset $Z$ of $X$ such that  $|\rho(x)|\leq\varepsilon$ for any $x\notin Z$. For any $z\in Z$ choose a finite  $\varepsilon$-net $H_z$ for the  interval $[-\rho(z),\rho(z)]$ and let
\[
H=\{g\in G:\ g(z)\in H_z\ \text{ for }\ z\in Z \text{ and }\ g(x)=0\ \text{ for }\ 
x\notin Z \}. 
\]
Then $H$ is a finite $\varepsilon$-net for $G$. Hence $G$ is a complete and totally
bounded metric space and so, $G$ is compact. Therefore $F=G+h$ is also compact. 
\quad $\Box$.

\section{Approximative solutions of difference equations} 

In this section we establish fundamental properties of the iterated remainder 
operator $r^m$ in the discrete case. Next, in Theorem \ref{T1}, we obtain our 
first main result. 

\begin{lemma}\label{dL2}
Assume $m\in\N$, $x\in\SQ$ and 
\begin{equation}\label{dEL1}
\sum_{n=1}^\infty n^{m-1}|x_n|<\infty. 
\end{equation}
Then there exists exactly one sequence $z\in\SQ$ such that 
\begin{equation}\label{Dmz}
z_n=\o(1) \qquad \text{and} \qquad \D^mz=x. 
\end{equation}
The sequence $z$ is defined by 
\begin{equation}\label{zn}
z_n=(-1)^m\sum_{j=n}^\infty\frac{(j-n+1)^{\overline{m-1}}}{(m-1)!}x_j.
\end{equation}
Moreover, if $k\in\N(0,m-1)$, then 
\begin{equation}\label{Dkzn}
\D^kz_n=(-1)^{m-k}\sum_{j=n}^\infty
\frac{(j-n+1)^{\overline{m-1-k}}}{(m-1-k)!}x_j.
\end{equation}
\end{lemma}
\textbf{Proof.} 
By Lemma 4 in \cite{J Migda 2014 b}, there exists exactly one sequence $z\in\SQ$ 
such that \eqref{Dmz} is satisfied. The sequence $z$ is defined by 
\begin{equation}\label{znb}
z_n=(-1)^m\sum_{k=0}^\infty\frac{(k+1)(k+2)\cdots(k+m-1)}{(m-1)!}x_{n+k}.
\end{equation}
Replacing $n+k$ by $j$ in \eqref{znb} we obtain \eqref{zn}. Moreover, \eqref{Dkzn} 
is a consequence of the proof of Lemma 3 in \cite{J Migda 2014 b}. \quad $\Box$.

\begin{lemma}\label{dL3}
Assume $m\in\N$, $\alpha\in(-\infty,0]$, $x\in\SQ$, 
\[
\sum_{n=1}^\infty n^{m-1-\alpha}|x_n|<\infty
\]
and $z\in\SQ$ is defined by \eqref{zn}. Then $z_n=\o(n^\alpha)$. 
\end{lemma}
\textbf{Proof.} 
The assertion follows from Lemma 4.2 of \cite{J Migda 2014 a}. \quad $\Box$.

\begin{definition}\label{def1}
Assume $m,p\in\N$, $x\in\F(\N_p)$ and 
\[ 
\sum_{n=p}^\infty n^{m-1}|x_n|<\infty. 
\]
We define $r^mx\in\F(\N_p)$ by 
\begin{equation}\label{rmxn}
r^mx_n=\sum_{j=n}^\infty\frac{(j-n+1)^{\overline{m-1}}}{(m-1)!}x_j.
\end{equation}
\end{definition}

\begin{lemma}\label{rmL}
Assume $m,p\in\N$, $k\in\N(0,m)$, $\alpha\in(-\infty,0]$,  
\[ 
\sum\limits_{n=p}^\infty n^{m-1-\alpha}|x_n|<\infty,  
\]
and $n\in\N_p$. Then we have 
\begin{equation}\label{rm|x|}
r^m|x|_n\leq\sum_{j=n}^\infty j^{m-1}|x_j|,
\end{equation}
\begin{equation}\label{Dkrm}
\D^kr^mx=(-1)^kr^{m-k}x=\o(n^{\alpha-k}),
\end{equation}
\begin{equation}\label{Dmrmx}
\D^mr^mx=(-1)^mx, \qquad r^mx=\o(n^{\alpha}).
\end{equation}
\end{lemma}
\textbf{Proof.} 
Obviously, \eqref{rm|x|} is a consequence of \eqref{rmxn}. If $k\in\N(0,m-1)$, 
then, by \eqref{zn} and \eqref{Dkzn}, we have 
\[
\D^kr^mx=(-1)^{-k}r^{m-k}x=(-1)^kr^{m-k}x. 
\]
Moreover, by \eqref{Dmz}, we get $\D^mr^mx=(-1)^mx$. Using Lemma \ref{dL3} and the  equality 
\[
m-1-\alpha=(m-k)-1-(\alpha-k)
\]
we have 
\[
r^{m-k}x=\o(n^{\alpha-k}). 
\]
Hence, we obtain \eqref{Dkrm}. Taking $k=m$ and $k=0$ in \eqref{Dkrm} 
we obtain \eqref{Dmrmx}. \quad $\Box$

\begin{theorem}\label{T1}
Assume $m,p\in\N$, $U\subset\R$, $M\geq 1$, $\mu>0$, 
$\alpha\in(-\infty,0]$, 
\begin{equation}\label{function f}
f:\N_p\times\R\to\R, \qquad \|f|\N_p\times U\|\leq M, 
\end{equation}
$f|\N_p\times U$ is continuous, $a,b,y\in\F(\N_p)$, $\D^my=b$, 
\begin{equation}\label{Mmu}
M\sum_{n=p}^\infty n^{m-1-\alpha}|a_n|\leq\mu, \quad \text{and} \quad 
y(\N_p)\subset\Int(U,\mu).
\end{equation}
Then there exists a sequence $x\in\F(\N_p)$ such that $x_n=y_n+\o(n^\alpha)$ and 
\[
\D^mx_n=a_nf(n,x_n)+b_n 
\]
for $n\geq p$.
\end{theorem}
\textbf{Proof.} 
Let $\rho=r^m|a|$ and 
\begin{equation}\label{Q}
Q=\{x\in\F(\N_p):\ |x-y|\leq M\rho\}. 
\end{equation}
Let $x\in Q$ and $n\in\N_p$. Then, using \eqref{rm|x|} and \eqref{Mmu}, we have 
\[
|x_n-y_n|\leq M\rho_n=Mr^m|a|_n\leq M\sum_{j=n}^\infty j^{m-1}|a_j|\leq\mu.
\]
Hence, using the inclusion $y(\N_p)\subset\Int(U,\mu)$, we obtain 
$x(\N_p)\subset U$. Therefore, by \eqref{function f},  
\begin{equation}\label{fM}
|f(n,x_n)|\leq M
\end{equation}
for any $x\in Q$ and $n\in\N_p$. For $n\geq p$ let 
\[
\bar{x}_n=a_nf(n,x_n). 
\] 
Then $|\bar{x}|\leq M|a|$. Thus $r^m|\bar{x}|\leq Mr^m|a|=M\rho$.  
Now, we define a sequence $Ax$ by 
\[
Ax=y+(-1)^mr^m\bar{x}. 
\]
Then 
\[
|Ax-y|=|r^m\bar{x}|\leq r^m|\bar{x}|\leq M\rho.
\]
Hence, by \eqref{Q}, $Ax\in Q$. Thus 
\[
AQ\subset Q. 
\]
Let $\varepsilon>0$. There exist $q\in\N_p$ and $\alpha>0$ such that 
\[
M\sum_{n=q}^\infty n^{m-1}|a_n|<\varepsilon \qquad\text{and}\qquad
\alpha\sum_{n=p}^{q}n^{m-1}|a_n|<\varepsilon.
\]
Let
\[
W=\{(n,t):\ n\in\N(p,q),\quad |t-y_n|\leq\mu\}
\]
Then $W$ is compact and, by \eqref{Mmu}, $W\subset\N_p\times U$. Hence $f|W$ is uniformly  continuous. Therefore there exists $\delta>0$ such that if $(n,t_1)$, $(n,t_2)\in W$ 
and $|t_1-t_2|<\delta$, then 
\[
|f(n,t_1)-f(n,t_2)|<\alpha. 
\] 
Choose $z\in Q$ such that\: $\|x-z\|<\delta$. Let 
$\bar{z}=a_nf(n,z_n)$. Then 
\[
\|Ax-Az\|=\sup_{n\geq p}|r^m\bar{x}_n-r^m\bar{z}_n|=
\sup_{n\geq p}|r^m(\bar{x}-\bar{z})_n|\leq r^m|\bar{x}-\bar{z}|_p
\]
\[
\leq\sum_{n=p}^\infty n^{m-1}|\bar{x}_n-\bar{z}_n|\leq
\sum_{n=p}^qn^{m-1}|\bar{x}_n-\bar{z}_n|+
\sum_{n=q}^\infty n^{m-1}|\bar{x}_n-\bar{z}_n|
\]
\[
\leq\alpha\sum_{n=p}^qn^{m-1}|a_n|+2M\sum_{n=q}^\infty n^{m-1}|a_n|<3\varepsilon.
\]
Hence, the mapping $A: Q\rightarrow Q$ is continuous. By Theorem \ref{rho}, 
$Q$ is a compact subset of the metric space 
\[
\F(\N_p)_y=y+\F_b(\N_p). 
\] 
Hence $AQ\subset Q$ is totally bounded. Therefore, by Theorem \ref{GST}, there exists 
$x\in Q$ such that $Ax=x$. Then 
\[
x=y+(-1)^mr^m\bar{x} 
\]
and, by \eqref{Dmrmx}, 
\[
\D^mx=\D^my+\D^m((-1)^mr^m\bar{x})=b+\bar{x}. 
\]
Therefore 
\[
\D^mx_n=a_nf(n,x_n)+b_n
\]
for $n\geq p$. Moreover, using \eqref{Dmrmx}, we have 
\[
x=y+(-1)^mr^m\bar{x}=y+\o(n^\alpha).
\]
The proof is complete. \quad $\Box$.

\section{Approximative solutions of differential equations}

In this section we establish fundamental properties of the iterated remainder 
operator $r^m$ in the continuous case. Next, in Theorem \ref{T2}, we obtain our 
second main result. 

\begin{lemma}\label{L1}
Let $a\in\R$, $b\in(a,\infty)$, $f\in\C[a,\infty)$, $m\in\N(0)$. Then 
\[
\int_a^b\int_t^b\frac{(s-t)^m}{m!}f(s)dsdt=
\int_a^b\frac{(s-a)^{m+1}}{(m+1)!}f(s)ds.
\]
\end{lemma}
\textbf{Proof.} 
Let $H:[a,b]\times[a,b]\to\R$, 
\[
H(t,s)=\frac{(s-t)^m}{m!}f(s).
\]
Then $H$ is continuous and 
\[
\int_a^b\int_t^b\frac{(s-t)^m}{m!}f(s)dsdt=\int_a^b\int_t^bH(t,s)dsdt
\]
\[
=\int_a^b\int_a^sH(t,s)dtds=\int_a^bf(s)\int_a^s\frac{(s-t)^m}{m!}dtds
\]
\[
=\int_a^bf(s)\left[-\frac{(s-t)^{m+1}}{(m+1)!}\right]_a^sds=
\int_a^b\frac{(s-a)^{m+1}}{(m+1)!}f(s)ds. 
\]
\quad $\Box$.

\begin{lemma}\label{L2}
Assume $m\in\N$, $t_0\in[0,\infty)$, $f\in\C[t_0,\infty)$ and 
\begin{equation}\label{EL1}
\int_{t_0}^\infty s^{m-1}|f(s)|ds<\infty. 
\end{equation}
Then there exists exactly one function $F:[t_0,\infty)\to\R$ such that 
\begin{equation}\label{Fmf}
F=\o(1) \qquad \text{and} \qquad F^{(m)}=f. 
\end{equation}
The function $F$ is defined by 
\begin{equation}\label{Ft}
F(t)=(-1)^m\int_t^\infty\frac{(s-t)^{m-1}}{(m-1)!}f(s)ds.
\end{equation}
Moreover, if $k\in\N(0,m-1)$, then 
\begin{equation}\label{Fkt}
F^{(k)}(t)=(-1)^{m-k}\int_t^\infty\frac{(s-t)^{m-1-k}}{(m-1-k)!}f(s)ds.
\end{equation}
\end{lemma}
\textbf{Proof.} 
Let $\varphi_0=\psi_0=f$. For $k\in\N(1,m)$ let 
$\varphi_k, \psi_k:[t_0,\infty)\to\R$, 
\[
\varphi_k(t)=\int_t^\infty\varphi_{k-1}(s)ds, \qquad 
\psi_k(t)=\int_t^\infty\frac{(s-t)^{k-1}}{(k-1)!}f(s)ds.
\]
By (\ref{EL1}), the integrals $\psi_k$ are convergent. Using Lemma \ref{L1} it 
is easy to see that $\varphi_k=\psi_k$ for $k\in\N(0,m)$. Let 
$F=(-1)^m\psi_m$. Obviously $\varphi_k'=-\varphi_{k-1}$ for $k\in\N(1,m)$. 
Hence 
\begin{equation}\label{Fk}
F^{(k)}=(-1)^{m-k}\psi_{m-k} 
\end{equation}
for $k\in\N(1,m)$. By (\ref{EL1}), we have 
\[
\int_{t_0}^\infty|\psi_{m-1}(s)|ds<\infty.
\]
Hence 
\[
\psi_m(t)=\int_t^\infty\psi_{m-1}(s)ds=\o(1). 
\]
Therefore $F=\o(1)$. Taking $k=m$ in \eqref{Fk} we obtain \eqref{Fmf}. Moreover, 
\eqref{Fkt} is a consequence of \eqref{Fk}. Now assume 
\[
G:[t_0,\infty)\to\R, \quad G^{(m)}=f \quad \text{and} \quad G=\o(1). 
\]
Then $(G-F)^{(m)}=0$ and so $G-F\in\Pol(m-1)$. Moreover $G-F=\o(1)$. 
Since $\Pol(m-1)\cap\o(1)=0$, we obtain $G-F=0$. 
The proof is complete. \quad $\Box$.

\begin{lemma}\label{L3}
Assume $m\in\N$, $t_0\in[0,\infty)$, $\alpha\in(-\infty,0]$, 
$f\in\C[t_0,\infty)$, 
\[
\int_{t_0}^\infty s^{m-1-\alpha}|f(s)|ds<\infty
\]
and $F:[t_0,\infty)\to\R$ is defined by \eqref{Ft}. Then $F(t)=\o(t^\alpha)$. 
\end{lemma}
\textbf{Proof.} 
Let $g,G:[t_0,\infty)\to\R$, 
\[
g(s)=s^{-\alpha}f(s), \qquad G(t)=\int_t^\infty\frac{(s-t)^{m-1}}{(m-1)!}|g(s)|ds.
\]
Then, by Lemma \ref{L2}, $G=\o(1)$. Moreover, 
\[
|t^{-\alpha}F(t)|=
t^{-\alpha}\left|\int_t^\infty\frac{(s-t)^{m-1}}{(m-1)!}f(s)ds\right|=
\left|\int_t^\infty\frac{(s-t)^{m-1}}{(m-1)!}t^{-\alpha}f(s)ds\right|
\]
\[
\leq\int_t^\infty\frac{(s-t)^{m-1}}{(m-1)!}t^{-\alpha}|f(s)|ds\leq 
\int_t^\infty\frac{(s-t)^{m-1}}{(m-1)!}s^{-\alpha}|f(s)|ds=G(t)=\o(1).
\]
Hence $F(t)=\o(t^\alpha)$. \quad $\Box$.

\begin{definition}\label{def2}
Assume $m\in\N$, $f\in\C[t_0,\infty)$ and 
\[
\int_{t_0}^\infty s^{m-1}|f(s)|ds<\infty. 
\]
We define $r^mf:[t_0,\infty)\to\R$ by
\begin{equation}\label{rmf-def}
(r^mf)(t)=\int_t^\infty\frac{(s-t)^{m-1}}{(m-1)!}f(s)ds.
\end{equation}
\end{definition} 

\begin{lemma}\label{L4} 
Assume $m\in\N$, $k\in\N(0,m)$, $\alpha\in(-\infty,0]$, $f\in\C[t_0,\infty)$ and 
\[
\int_{t_0}^\infty s^{m-1-\alpha}|f(s)|ds<\infty. 
\]
Then for $t\geq t_0$ we have 
\begin{equation}\label{rm|f|}
r^m|f|(t)\leq\int_{t}^\infty s^{m-1}|f(s)|ds 
\end{equation}
Moreover, 
\begin{equation}\label{rmfk}
(r^mf)^{(k)}=(-1)^kr^{m-k}f=\o(t^{\alpha-k}),
\end{equation}
\begin{equation}\label{rmfm}
(r^mf)^{(m)}=(-1)^mf, \qquad r^mf=\o(t^{\alpha}).
\end{equation}
\end{lemma}
\textbf{Proof.} 
Inequality \eqref{rm|f|} is an easy consequence of \eqref{rmf-def}. 
By \eqref{rmf-def}, \eqref{Ft} and \eqref{Fkt}, we have 
\[
(r^mf)^{(k)}=(-1)^kr^{m-k}f.
\]
Using Lemma \ref{L3}, we obtain the right equality in \eqref{rmfk}. 
Taking $k=m$ and $k=0$ in \eqref{rmfk} we obtain \eqref{rmfm}.
\quad $\Box$.

\begin{theorem}\label{T2}
Assume $m\in\N$, $U\subset\R$, $M\geq 1$, $t_0\in[0,\infty)$, $\mu>0$, 
$\alpha\in(-\infty,0]$, 
\begin{equation}\label{cfunction f}
I=[t_0,\infty), \qquad f:I\times\R\to\R, \qquad \|f|I\times U\|\leq M,
\end{equation}
$f|I\times U$ is continuous, $a,b\in\C(I)$, $y\in\C^m(I)$, $y^{(m)}=b$, 
\begin{equation}\label{int}
M\int_{t_0}^\infty s^{m-1-\alpha}|a(s)|ds\leq\mu, \quad \text{and} \quad 
y(I)\subset\Int(U,\mu).
\end{equation}
Then there exists a function $x\in\C^m(I)$ such that $x(t)=y(t)+\o(t^\alpha)$ and 
\[
x^{(m)}(t)=a(t)f(t,x(t))+b(t) 
\]
for $t>t_0$.
\end{theorem}
\textbf{Proof.} 
Define a function $\rho$ and a subset $Q$ of $\C(I)$ by 
\[
\rho=r^m|a|, \qquad Q=\{x\in\C(I):\ |x-y|\leq M\rho\}. 
\] 
By \eqref{rmf-def} and \eqref{int} we have $M\rho\leq\mu$. 
Let $x\in Q$. Then $|x-y|\leq\mu$. 
Hence, using the inclusion $y(I)\subset\Int(U,\mu)$, we get $x(I)\subset U$. 
Therefore, by \eqref{cfunction f}, we have 
\begin{equation}\label{ftxM}
|f(t,x(t))|\leq M
\end{equation}
for any $x\in Q$ and $t\in I$. For $t\geq t_0$ let 
\begin{equation}\label{barx}
\bar{x}(t)=a(t)f(t,x(t)). 
\end{equation}
Then $|\bar{x}|\leq M|a|$. Thus $r^m|\bar{x}|\leq Mr^m|a|=M\rho$.  
Now, we define a function $Ax$ by 
\[
Ax=y+(-1)^mr^m\bar{x}. 
\]
Then 
\[
|Ax-y|=|r^m\bar{x}|\leq r^m|\bar{x}|\leq M\rho.
\]
Hence $Ax\in Q$. Thus 
\[
AQ\subset Q. 
\]
Let $\varepsilon>0$. By \eqref{int}, there exist $k\in I$ and $\alpha>0$ such that 
\[
M\int_{k}^\infty s^{m-1}|a(s)|ds<\varepsilon \qquad\text{and}\qquad
\alpha\int_{t_0}^{k}s^{m-1}|a(s)|ds<\varepsilon.
\]
Let
\[
W=\{(t,s):\ t\in[t_0,k],\quad |s-y(t)|\leq\mu\}
\]
Then $W$ is compact and $f|W$ is uniformly continuous. Hence there 
exists $\delta>0$ such that if $(t,s_1)$, $(t,s_2)\in W$ and 
$|s_1-s_2|<\delta$, then 
\[
|f(t,s_1)-f(t,s_2)|<\alpha. 
\] 
Choose $z\in Q$ such that\: $\|x-z\|<\delta$. Let $\bar{z}=a(t)f(t,z(t))$. 
Then, using \eqref{rm|f|} and \eqref{ftxM}, we obtain 
\[
\|Ax-Az\|=\sup_{t\geq t_0}|(r^m\bar{x})(t)-(r^m\bar{z})(t)|=
\sup_{t\geq t_0}|r^m(\bar{x}-\bar{z})(t)|\leq r^m|\bar{x}-\bar{z}|(t_0)
\]
\[
\leq\int_{t_0}^\infty s^{m-1}|\bar{x}-\bar{z}|(s)ds\leq
\int_{t_0}^ks^{m-1}|\bar{x}-\bar{z}|(s)ds+
\int_k^\infty s^{m-1}|\bar{x}-\bar{z}|(s)ds
\]
\[
\leq\alpha\int_{t_0}^ks^{m-1}|a(s)|ds+2M\int_k^\infty s^{m-1}|a(s)|ds<3\varepsilon.
\]
Hence, the mapping $A: Q\rightarrow Q$ is continuous. 
Assume $x\in Q$. If $t_1,t_2\in I$, then 
\[
|A(x)(t_1)-A(x)(t_2)|=|r^m\bar{x}(t_1)-r^m\bar{x}(t_2)|
\]
\[
=\left|\int_{t_1}^\infty r^{m-1}\bar{x}(s)ds-
\int_{t_2}^\infty r^{m-1}\bar{x}(s)ds\right|=
\left|\int_{t_1}^{t_2}r^{m-1}\bar{x}(s)ds\right|\leq
\int_{t_1}^{t_2}r^{m-1}|\bar{x}|(s)ds.
\]
Moreover, if $t\in I$, then 
\[
r^{m-1}|\bar{x}|(t)\leq\int_t^\infty s^{m-2}|\bar{x}(s)|ds\leq
\int_{t_0}^\infty s^{m-2}|a(s)f(s,x(s))|ds
\]
\[
\leq M\int_{t_0}^\infty s^{m-1}|a(s)|ds\leq\mu.
\]
Hence 
\[ 
|A(x)(t_1)-A(x)(t_2)|\leq\mu\int_{t_1}^{t_2}ds=|t_2-t_1|\mu. 
\]
Therefore the family $AQ$ is equicontinuous. By Theorem \ref{rho}, $Q$ is closed, 
convex, pointwise bounded, and stable at infinity. Hence $AQ$ is pointwise bounded,  equicontinuous and stable at infinity. By Theorem \ref{GAT}, $AQ$ is totally bounded. 
By Theorem \ref{GST}, there exists $x\in Q$ such that $Ax=x$. Then 
\begin{equation}\label{y+rmx}
x=y+(-1)^mr^m\bar{x}. 
\end{equation}
Using \eqref{rmfm}, we obtain 
\[
x^{(m)}=y^{(m)}+((-1)^mr^m\bar{x})^{(m)}=b+\bar{x}.
\]
Therefore, by \eqref{barx}, we get 
\[
x^{(m)}(t)=a(t)f(t,x(t))+b(t)
\]
for $t>t_0$. Moreover, using \eqref{y+rmx} and \eqref{rmfm}, we have 
\[
x=y+(-1)^mr^m\bar{x}=y+\o(t^\alpha).
\]
The proof is complete. \quad $\Box$.

\section{Remarks}

In this section, we give some additional remarks on the regional topology. In particular, 
in Remark \ref{top 03} we show that if $X$ is an extraordinary regional space, then 
the regional topology is almost linear but not linear. In Example \ref{ord ess} we 
show, that in Theorem \ref{GST}, the assumption of ordinarity of the set $AQ$ is 
essential. At the end of the section we show that if $X$ is a locally compact but 
noncompact metric space, then the notion of equicontinuity at infinity of a family 
$F\subset\F(X)$ is equivalent to the known notion of equiconvergence at infinity.
\smallskip\\
Let $X$ be a regional normed space. The following remark is a consequence of the 
fact, that $X$ is a topological disjoint union of all its regions.

\begin{remark}\label{top 01} 
A subset $Y$ of $X$ is closed in $X$ if and only if, $Y\cap X_p$ is closed in 
$X_p$ for any region $X_p$. If $Z$ is a topological space, then a map 
$\varphi:X\to Z$ is continuous if and only if $\varphi|X_p$ is continuous for any 
region $X_p$. If $Y\subset X$ and $p\in X$, then the set $Y\cap X_p$ is closed and open 
in $Y$.
\end{remark}
Let $Z$ be a linear subspace of $X$ such that 
\[
Z\cap X_0=\{0\},
\]
$q\in X$, and $Y=q+Z$. If $p\in Y$, then $Y\cap X_p=\{p\}$. Hence, by 
Remark \ref{top 01}, we obtain:

\begin{remark}\label{top 02} 
If $Z$ is a linear subspace of $X$ such that $Z\cap X_0=\{0\}$, $q\in X$, then the 
topology on $q+Z$ induced from $X$ is discrete.
\end{remark}

\begin{example}
A formula  
\[
\|(x,y)\|= 
\begin{cases}
|y| &\text{if $x=0$}\\
\infty &\text{if $x\neq 0$}
\end{cases}
\]
define a regional norm on the space $X=\R^2$. A subset $L$ of $X$ is a region in $X$ 
if and only if $L$ is a vertical line. If $L$ is a vertical line, then the 
topology induced on $L$ from $X$ is the usual Euclidean topology. On the other hand, 
the topology induced from $X$ on any nonvertical line $L$ is discrete. Moreover, if 
$\varphi:\R\to\R$ is an arbitrary function, then the topology induced on the graph 
\[
G(\varphi)=\{(x,\varphi(x)):\ x\in\R\}
\]
is discrete.
\end{example}

\begin{example}\label{[x,y]}
Let $x,y\in X$ and let $I$ denote the line segment 
\[
I=[x,y]=x+[0,y-x]=x+[0,1](y-x)=\{x+\lambda(y-x):\ \lambda\in[0,1]\}.
\]
If $\|y-x\|<\infty$, then 
\[
[0,y-x]=[0,1](y-x)\subset X_0. 
\]
Hence, in this case, $I$ is topologically equivalent to standard interval $[0,1]$. 
Moreover, $I$ is a closed subset of $X_x$ and, by Remark \ref{top 01}, $I$ is 
closed in $X$.
\smallskip\\
Assume $\|y-x\|=\infty$. Then 
\[
X_0\cap\R(y-x)=\{0\}
\]
and, by Remark \ref{top 02}, the topology induced on $x+\R(y-x)$ from $X$ is 
discrete. Moreover 
\[
I\subset x+\R(y-x).
\]
Hence, in this case, the topology induced on $I$ is discrete. If $p\in X$, then 
$I\cap X_p$ is one point or empty subset of $X_p$ and, by Remark \ref{top 01}, $I$ is 
closed in $X$.
\end{example}
Using Remark \ref{top 01} it is not difficult to see that the addition 
$X\times X\to X$ is continuous. Moreover, if $p\in X$, and $\lambda\in\R$, then the  translation and $x\mapsto p+x$ and the homothety $x\mapsto \lambda x$ are continuous 
maps from $X$ to $X$. Let $\mu:\R\times X\to X$ denote the multiplication by scalars.  
If there exists an $x\in X$ such that $\|x\|=\infty$, then $\mu(\R\times\{x\})=\R x$ 
and, by Remark \ref{top 02}, the topology induced on $\R x$ is discrete. Hence, the  restriction $\mu|\R\times\{x_0\}$ is discontinuous. 
This implies the discontinuity of $\mu$. Hence we have

\begin{remark}\label{top 03} 
The addition
\[
X\times X\to X, \qquad (x,y)\mapsto x+y
\]
is continuous. If $p\in X$, then the translation 
\[
X\to X, \qquad x\mapsto p+x 
\]
is a homeomorphism. If $\lambda$ is a nonzero scalar, the the homothety 
\[
X\to X, \qquad x\mapsto \lambda x 
\]
is a homeomorphism. If $X$ is extraordinary, then the multiplication by scalars 
\[
\R\times X\to X, \qquad (\lambda,x)\mapsto\lambda x
\]
is discontinuous.
\end{remark}
Every region $X_p$ is topologically equivalent to the normed space $X_0$. Hence $X_p$ is  connected. Moreover, $X_p$ is open and closed subset of $X$. Hence $X_p$ is a maximal  connected subset of $X$. Therefore any connected subset of $X$ is contained in certain  region. 
\smallskip\\ 
Obviously, $X$ is a Hausdorff space and so, any compact subset of $X$ is closed. 
Moreover, if $C\subset X$ is compact and $p\in X$, then $C_p=C\cap X_p$ is closed and 
open in $C$. Hence $C_p$ is compact and 
\[
C\subset C_{p_1}\cup\cdots\cup C_{p_n}
\]
for some $p_1,\dots,p_n\in X$.
\smallskip\\ 
Now, assume that $C\subset X$ is compact and convex. Let $x,y\in C$ and let $I=[x,y]$. 
If $\|y-x\|=\infty$, then, by Example \ref{[x,y]}, $I$ is closed and discrete subset 
of $C$. Hence $I$ is compact and discrete. It is impossible. This means that $C$ 
is ordinary. Obviously, any convex and ordinary subset of $X$ is connected.
\smallskip\\ 
Therefore we obtain:

\begin{remark}\label{top 04} 
Any connected subset of $X$ is ordinary. Any compact subset of $X$ is a finite sum 
of ordinary compact subsets. Any compact and convex subset of $X$ is connected and 
ordinary.
\end{remark}

\begin{example}\label{ord ess}
Let $X$ be a regional normed space, $x,y\in X$, $\|y-x\|=\infty$, 
\[
Q=[x,y], \qquad 
A:Q\to Q, \qquad Az= 
\begin{cases}
x &\text{for $z\in(x,y]$}\\
y &\text{for $z=x$}
\end{cases}.
\]
Then $Q$ is closed and convex, $A$ is continuous, $AQ=\{x,y\}$ is totally bounded and 
$Az\neq z$ for any $z\in Q$.
\end{example}
Assume that $X$ is a locally compact, noncompact metric space. 
\smallskip\\
We say that a real number $p$ is a limit at infinity of a function $f:X\to\R$ 
if for any $\varepsilon>0$ there exists a compact subset $Z$ of $X$ such that  $|f(t)-p|<\varepsilon$ for any $t\notin Z$. 
Then we write 
\[
p=\lim_{t\to\infty}f(t)
\]
and say that $f$ is convergent at infinity. 
We say that a family $F\subset\F(X)$ is equiconvergent at infinity if all 
functions $f\in F$ are convergent at infinity and for any $\varepsilon>0$ there 
exists a compact $Z\subset X$ such that  
\[
\left|f(s)-\lim_{t\to\infty}f(t)\right|<\varepsilon 
\]
for any $f\in F$ and $s\notin Z$. 
\smallskip\\
In the next remark we show that a family $F\subset\F(X)$ is equiconvergent 
at infinity if and only if, it is equcontinuous at infinity.

\begin{remark}\label{equiconvergence}
Obviously every equiconvergent at infinity family $F\subset F(X)$ is equcontinuous 
at infinity. If a family $F\subset F(X)$ is equicontinuous at infinity, then  
for any natural $n$ there exists a compact $Z_n\subset X$ such that $|f(s)-f(t)|<1/n$ 
for every $f\in F$ and $s,t\notin Z_n$. For $n\in\N$ let 
\[
K_n=Z_1\cup\dots\cup Z_n.
\]
Then $K_n$ is compact. 
Choose a sequence $(t_n)$ in $X$ such that $t_n\notin K_n$ for any $n$. 
If $f\in F$, then $(f(t_n))$ is a Cauchy sequence and there exists a limit 
$p$ of this sequence. It is easy to see that  
\[
p=\lim_{t\to\infty}f(t)
\]
and $F$ is equiconvergent at infinity.
\end{remark}

Note that if $X$ is a compact space, then every family $F\subset\F(X)$ is 
equicontinuous at infinity but, in this case, equiconvergence at infinity is not 
defined.

\end{document}